\documentclass{mcom-l}
\usepackage{amssymb}
\usepackage{url}

\newcommand{\Z}{\mathbb{Z}}
\newcommand{\C}{\mathbb{C}}
\newcommand{\F}{\mathbb{F}}
\newcommand{\Fc}{\mathcal{F}}
\newcommand{\Oc}{\mathcal{O}}
\newcommand{\HH}{\mathbb{H}}
\newcommand{\w}{\mathfrak{w}}
\newcommand{\wl}{\w_\ell}
\newcommand{\wt}{\w_2}
\newcommand{\f}{\mathfrak{f}}
\newcommand{\jbar}{\overline j}
\newcommand{\lf}{\mathfrak{l}}
\newcommand{\llaur}{((}
\newcommand{\rlaur}{))}
\newcommand{\divides}{|}
\newcommand{\legendre}[2]{\genfrac {(}{)}{1pt}{}{#1}{#2}}

\newtheorem {result}{Result}[section]
\theoremstyle {remark}
\newtheorem* {remark}{Remark}

\copyrightinfo {2008}{Andreas Enge}

\numberwithin{equation}{section}

\title{Computing modular polynomials in quasi-linear time}

\author {Andreas Enge}
\address {INRIA Saclay--\^Ile-de-France
\& Laboratoire d'In\-for\-ma\-tique (CNRS/UMR 7161),
\'Ecole polytechnique, 91128 Palaiseau Cedex, France}
\email {enge@lix.polytechnique.fr}


\date{April 24, 2007; revised May 5, 2008}

\begin{document}

\begin{abstract}
We analyse and compare the complexity of several algorithms for computing
modular polynomials. Under the assumption that rounding errors do not
influence the correctness of the result, which appears to be
satisfied in practice, we show that an algorithm relying on floating point
evaluation of modular functions and on interpolation has a complexity that
is up to logarithmic factors linear in the size of the computed polynomials.
In particular, it obtains the classical modular polynomial $\Phi_\ell$
of prime level~$\ell$ in time
\[
O \left( \ell^2 \log^3 \ell \, M (\ell) \right)
\subseteq O \left( \ell^3 \log^{4 + \varepsilon} \ell \right),
\]
where $M (\ell)$ is the time needed to multiply two $\ell$-bit
numbers.

Besides treating modular polynomials for $\Gamma^0 (\ell)$, which are an
important ingredient in many algorithms dealing with isogenies of elliptic
curves, the algorithm is easily adapted to more general situations. Composite
levels are handled just as easily as prime levels, as well as polynomials
between a modular function and its transform of prime level, such as the
Schl\"afli polynomials and their generalisations.

Our distributed implementation of the algorithm confirms the theoretical
analysis by computing modular equations of record level around $10000$ in less
than two weeks on ten processors.
\end{abstract}

\maketitle

\section{Definitions and main result}

Modular polynomials, in their broadest sense, are bivariate polynomials
with a pair of modular functions as zero. Given any two modular
functions $f$ and~$g$ for arbitrary congruence subgroups, the function
fields $\C (f)$ and $\C (g)$ are finite extensions of $\C (j)$, so that
there are two polynomials relating $f$ resp. $g$ to $j$. Taking the
resultant of these polynomials with respect to $j$, one sees that a
polynomial relationship between $f$ and $g$ exists. In practice, one is
rather interested in the minimal polynomial of $f$ over $\C (g)$, say,
that will be called the modular polynomial of $f$ with respect to $g$.
If the functions satisfy some conditions on the rationality and
integrality of their $q$-expansion coefficients, the
modular polynomial has rational integral coefficients.

Different modular polynomials parameterise moduli spaces related to
elliptic curves. Let $\Gamma = \operatorname{Sl_2} (\Z) / \{ \pm 1 \}
= \operatorname{PSl_2} (\Z)$ be the full modular group, and $\C_\Gamma =
\C (j)$ the field of modular functions invariant under $\Gamma$; $j$
itself parameterises isomorphism classes of elliptic curves. Of special
interest for applications is the congruence subgroup
\[
\Gamma^0 (\ell) =
\begin{pmatrix} 1 & 0 \\ 0 & \ell \end{pmatrix}^{-1} \Gamma
\begin{pmatrix} 1 & 0 \\ 0 & \ell \end{pmatrix} \cap \Gamma
= \left\{ \begin{pmatrix} a & b \\ c & d \end{pmatrix} \in \Gamma : \ell
\divides b \right\}
\]
for $\ell$ prime; its modular polynomial parameterises isomorphism
classes of elliptic curves together with an isogeny of degree $\ell$.
This is heavily used in Atkin's and Elkies's improvements to Schoof's
algorithm for counting points on elliptic curves \cite{elk98}.
Other applications include the computation of the endomorphism ring of
an elliptic curve \cite{koh96,fm02} or of an isogeny between two
elliptic curves \cite{gal99}.

More general modular polynomials occur in modern complex multiplication
constructions for elliptic curves \cite{es04,em07,em02}; they are
discussed in Section~\ref{sec:further}. When counting points on an
elliptic curve by the Schoof--Elkies--Atkin algorithm, the modular
polynomials are usually supposed to be available after a precomputation
phase. This makes sense in situations where the algorithm is carried out
many times for curves of about the same size, as is typical in
cryptography. When establishing point counting records
\cite{em06}, however, one quickly realises that the computation of the
modular polynomials itself becomes a bottleneck that may even dominate
from a complexity theoretic point of view.

In this article, after surveying in Section~\ref{sec:approaches} the
state of the art as it appears in the literature, I present an algorithm
based on evaluation and interpolation that computes modular polynomials
in time essentially linear in the bit size of the polynomial. As it is
based on floating point computations, one has to assume that rounding
errors do not disturb the output, a heuristic that is supported by the
implementation described in Section~\ref{sec:implementation}. A precise
analysis of rounding errors to obtain a rigorously proved result appears to
be out of reach. Note, however, that the correctness of the output may be
checked probabilistically. For instance, one may instantiate the variable $j$
occurring in the modular polynomial by the j-invariant of an elliptic curve
over a finite field and verify that the specialised polynomial has the
expected splitting behaviour.

Under the assumption that it is sufficient to use a floating point precision
of $O (n)$ for a polynomial whose largest coefficient has $n$ digits, the
following holds:
\begin{result}[heuristic]
\label {th}
Let $\Gamma' \subseteq \Gamma$ be a congruence subgroup, $f$ a
modular function for $\Gamma'$ and $\Phi (X) \in \C (j)[X]$
the characteristic polynomial of $f$ with respect to the function
field extension $\C_{\Gamma'} / \C (j)$. Assume that $\Phi$ has coefficients
in $\Z [j]$, and that the largest integer coefficient occurring in it
has $n$ digits.
Suppose that a system of representatives of
$\Gamma' \backslash \Gamma$ is known and that $f$ can be evaluated at
precision $O (n)$ in time $O (\log n \, M (n))$,
where $M (n) \subseteq O \left( n \log^{1 + \varepsilon} \right)$
is the time needed to multiply two numbers with $n$ digits.
Then there is an algorithm that computes $\Phi$ in time
\[
O \left( \deg_X \Phi \, \deg_j \Phi \,
\big( \log^2 \max (\deg_X \Phi, \deg_j \Phi) + \log n \big)
\, M (n) \right).
\]
In particular, the classical modular polynomial $\Phi_\ell$ such that
$\Phi_\ell (j (z), j (\ell z)) = 0$ for $\ell$ prime is obtained in time
\[
O \left( \ell^2 \log^2 \ell \, M (\ell \log \ell) \right)
\subseteq O \left( \ell^3 \log^{4 + \varepsilon} \ell \right).
\]
\end {result}

If $\Phi$ is a dense polynomial and all of its coefficients (or at least
a constant fraction of them) are of bit size about $n$, then this
algorithm is linear in the size of the polynomial, except for the
logarithmic factors.

The time bound of $O (\log n \, M (n))$ for evaluating modular functions
at precision~$n$ is motivated by an algorithm due to Dupont;
it relies on Newton iterations on an expression involving the
arithmetic-geometric mean and reaches this complexity for a wide class
of modular functions, including those built from Dedekind's
$\eta$ function and in particular $k$, $k'$ and $j$
\cite[Theorem~5 and Section~7.2]{dup07}.
Alternatively, one may use an algorithm relying on multievaluation of
$q$-expansions as described in \cite[Section~6.3]{eng06}. Its
complexity, worse by a factor of $\log n$, is still
sufficient for the running time of Result~\ref {th} to hold.
The bound $M (n) \in O \left( n \log^{1 + \varepsilon} \right)$ may
be achieved by the classical Sch\"onhage--Strassen algorithm \cite{ss71}
or the more recent and asymptotically faster algorithm due to
F\"urer \cite{fue07}.

The roots of $\Phi$, that is, the algebraic conjugates of $f$, are given
by the $f (M_\nu z)$ with $\Gamma' M_\nu$ running through a system of
representatives of the residue classes $\Gamma' \backslash \Gamma$
(see \cite{deu58}), so that
\begin {equation}
\label {eq:haupt}
\Phi (X) = \prod_{\nu = 1}^{[\Gamma : \Gamma']} (X - f (M_\nu z))
= X^{[\Gamma : \Gamma']} + \sum_{r=1}^{[\Gamma : \Gamma']}
c_r X^{[\Gamma : \Gamma'] - r},
\end {equation}
where the coefficients $c_r$ of $\Phi$ are the elementary symmetric
functions of the conjugates $f (M_\nu z)$.

By the $q$-expansion principle (cf. \cite[\S3]{deu58}), a sufficient
condition for the $c_r$ to be polynomials in $\Z [j]$ is that
$f$ has rational integral $q$-coefficients, is
holomorphic in $\HH = \{ z \in \C : \Im z > 0 \}$ and all conjugates
$f \left( \frac {az + b}{c z + d} \right)$ with $\begin{pmatrix} a & b
\\ c & d \end{pmatrix} \in \Gamma$ have integral algebraic
$q$-coefficients.

\section{Approaches to computing modular polynomials for $\Gamma^0
(\ell)$}
\label{sec:approaches}

The case of $\Gamma^0 (\ell)$ with $\ell$ prime, which is the most
important one for applications, is also the simplest one from a
theoretical point of view: The residue classes $\Gamma^0 (\ell)
\backslash \Gamma$ are represented by
\[
T^\nu = \begin{pmatrix} 1 & \nu \\ 0 & 1 \end{pmatrix},
\nu = 0, \ldots, \ell - 1,
\text{ and }
S = \begin{pmatrix} 0 & -1 \\ 1 & 0 \end{pmatrix}.
\]
(The literature often concentrates on
$\Gamma_0 (\ell) = \begin{pmatrix} \ell & 0 \\ 0 & 1 \end{pmatrix}^{-1}
\Gamma \begin{pmatrix} \ell & 0 \\ 0 & 1 \end{pmatrix} \cap \Gamma$.
Since $\Gamma_0 (\ell) = S^{-1} \Gamma^0 (\ell) S$, a modular function
$f$ for $\Gamma^0 (\ell)$ yields the function $f (Sz)$ for $\Gamma_0
(\ell)$, and $f (Sz)$ is actually a conjugate of $f$. So both result in
the same modular polynomial $\Phi$. Or phrased differently:
$\C_{\Gamma^0 (\ell)} / \C_\Gamma$ is not normal, and another one of its
embeddings is $\C_{\Gamma_0 (\ell)} / \C_\Gamma$.)

Different functions have been suggested in the literature:
\begin{itemize}
\item
$j (z/\ell)$, leading to the so-called ``classical'' or ``traditional''
modular polynomials;
\item
$\wl (z)^{2s} = \left( \frac {\eta (z / \ell)}{\eta (z)} \right)^{2s}$ for
$s = 12 / \gcd (12, \ell - 1)$, yielding the so-called ``canonical''
modular polynomials;
\item
functions that are moreover invariant under the Fricke--Atkin--Lehner
involution
$z \mapsto \frac {- \ell}{z}$; in particular, functions of the form
\[
\frac {T_r  (\eta (z) \eta (\ell z))}{\eta (z) \eta (\ell z)}
\]
evaluated in $\frac {-1}{z}$ with $T_r$ the $r$-th Hecke operator and
linear combinations thereof as described in \cite[Ch.~5.3]{mul95alg}.
\end{itemize}

All these functions lead to modular polynomials with coefficients in
$\Z$.

\subsection {Heights of modular polynomials}
\label {ssec:height}

The (logarithmic) height of a modular polynomial with coefficients in $\Z$ is
defined as the logarithm of the largest absolute value of the coefficients. It
provides a lower bound on the arithmetic precision required to compute the
polynomial. For the classical polynomial between $j (z)$ and $j (z / \ell)$, the
height is given by
\[
6 (\ell + 1) \left( \log \ell + O (1) \right)
\subseteq O (\ell \log \ell)
\]
according to \cite {coh84}. For alternative modular polynomials, one observes a
similar growth of the coefficients with $\ell$, but the constant is usually
smaller. It should be possible to adapt the argumentation in \cite {coh84} to
obtain a similar bound for further classes of modular polynomials.

\subsection {Using $q$-expansions}
\label {ssec:qexpansions}

The system of representatives of $\Gamma^0 (\ell) \backslash \Gamma$ and thus
the conjugates of the modular function $f$ being explicitly known, one may use
\eqref {eq:haupt} to compute the modular polynomial $\Phi$, if
only a procedure for recognising its coefficients $c_r$ as polynomials in $\Z
[j]$ is available. The straightforward and most popular approach is to use
$q$-expansions: The expansion of $j$ is of the form $q^{-1} + \cdots$,
so that knowing the expansion $c_r = c_{r,k} q^{-k} + \cdots$,
one deduces that $c_r$ has leading term $c_{r, k} j^k$ as a polynomial in $j$;
one subtracts the $q$-expansion of this term and continues in the same vein.
In fact, only the non-positive $q$-powers in the series expansions
of the $c_r$ are needed.
Several variants of the algorithm are described in \cite {mor95} with a brief
complexity analysis in \cite [Section 3]{elk98}. We shall describe only the
algorithm with the best complexity and give a detailed analysis of its running
time.

Let $f$ be a modular function for $\Gamma^0 (\ell)$ with valuation
$v < 0$ at infinity; that is, $f$ admits a $q$-expansion of the form
\[
f (z) = \sum_{k = v \ell}^\infty a_k q^{k/\ell} \in \C \llaur q^{1/\ell}
\rlaur
\text { with }
q^{1/\ell} = e^{2 \pi i z/\ell}.
\]
Suppose that the $a_k$ are rational integers. (For instance, $f = j (z/\ell)$
resp. $f = \wl^{2s}$ are valid choices with $v = -1/\ell$ resp.
$v = - \frac {s (\ell - 1)}{12 \ell}$.)
Then the $q$-expansions of the conjugates $f_\nu = f (T^\nu z)$ are of
the same form with the coefficients $a_k$ being multiplied by powers of
$\zeta_\ell = e^{2 \pi i / \ell}$.

The last conjugate $f_\infty = f (S z)$ poses difficulties since it
cannot be described generically, $q \left( \frac {-1}{z} \right) = q^{-
2 \pi i / z}$ being unrelated to $q$ in a simple way, and is thus
treated separately. (And it is the fact that there is only one of them that
makes the case of $\Gamma' = \Gamma^0 (\ell)$ stand out so that it is
comparatively easily handled via the $q$-expansion approach.)
Denote by $v_\infty$ the order of $f_\infty$ at infinity.
For $f = j (z / \ell)$, one has $f_\infty = j (\ell z) = q^{- \ell} + \cdots$
and $v_\infty = - \ell$, and for $f = \wl^{2s}$, the transformation behaviour
of $\eta$ yields $f_\infty = \left( \sqrt \ell \, \frac {\eta (\ell z)}{\eta (z)}
\right)^{2s} = \ell^s q^{s (\ell - 1) / 12} + \cdots$ with
$v_\infty = \frac {s (\ell - 1)}{12} = \frac {\ell - 1}{\gcd (12, \ell -1)}$.

In particular, $v_\infty$ may be positive or negative. If it is negative, the
term with lowest $q$-exponent occurring in \eqref {eq:haupt} is
$f_\infty \cdot \prod_{\nu = 0}^{\ell - 1} f_\nu$ of exponent $\ell v +
v_\infty$; so the degree of $\Phi$ in $j$ is $|\ell v + v_\infty|$, reached
for $c_{\ell + 1}$. For instance, the degree in $j$ is $\ell + 1$ for the
classical modular polynomial $\Phi_\ell$. If $v_\infty$ is positive, the term
with lowest $q$-exponent is $\prod_{\nu = 0}^{\ell - 1} f_\nu$, so the degree
$|\ell v|$ of $\Phi$ in $j$ is reached in $c_\ell$. For instance, the degree
in $j$ of the canonical modular polynomials is $\frac {s (\ell - 1)}{12}$.

We first consider the complexity of determining the partial product
\[
\prod_{\nu = 0}^{\ell - 1} (X - f_\nu) = X^\ell + \sum_{r=1}^\ell \tilde c_r
X^{\ell - r} \in \Z [\zeta_\ell] \llaur q^{1/\ell} \rlaur [X].
\]
One may proceed by computing the Newton sums
$s_r = \sum_{\nu = 0}^{\ell - 1} f_\nu^r$ for $1 \leq r \leq \ell$ and
use Newton's recurrence formul{\ae} to deduce the $\tilde c_r$. Following
\cite {mor95}, one obtains with $f^r = \sum_{k = r v \ell}^\infty a_{k,r}
q^{k/\ell}$ that
\begin {equation}
\label {eq:decimation}
s_r = \sum_{k = \lceil rv \rceil}^\infty
\ell \, a_{\ell k, r} \, q^k \in \Z \llaur q \rlaur,
\end {equation}
so that in fact the roots of unity do not interfere with the
computations. One sees that the valuations of the $s_r$ at infinity are
bounded below by $\lceil r v \rceil$ (and in general, they will be equal
to it). Assume that $d+1$ terms of them are known, that is, $s_r$ is known up
to and including the coefficient of $q^{\lceil r v \rceil + d}$. Then one
shows by induction on Newton's relation
\begin {equation}
\label {eq:newton}
\tilde c_r = \frac {1}{r} \sum_{k = 1}^r (-1)^{k+1} s_k \tilde c_{r-k}
\end {equation}
that the valuation of $\tilde c_r$ is also bounded below by $\lceil r v
\rceil$ and that it is also known up to $q^{\lceil r v \rceil + d}$.

If $v_\infty \geq 0$, we need only the terms of the $\tilde c_r$ with
non-positive exponents; for $v_\infty < 0$, we loose some precision,
in fact $|v_\infty|$ terms, through the final multiplication by $X -
f_\infty$. Thus, the value required for $d$ is given by
\[
d = \ell |v| + \max (0, -v_\infty) = \deg_j \Phi.
\]
This is an indication that the algorithm given above should be optimal among
those relying on $q$-expansions: Outputting polynomials with up to $\deg_j
\Phi + 1$ coefficients, it manipulates series with as many terms.
Unfortunately, the computation of the $s_r$ involves decimating the series for
$f^r$ by $\ell$. Eventually, a factor of $\ell$ is lost in the running time
since $O (\ell d)$ terms of the $f^r$ are needed.

Once the $q$-coefficients of $f$ are known, the $f^r$ and $s_r$ are thus
computed by $O (\ell)$ multiplications of series with $O (\ell d)$ terms and
the $\tilde c_r$ by $O (\ell^2)$ multiplications of series with $O (d)$ terms.
Let $M_q (d)$ be the number of
arithmetic operations in $\Z$ required to multiply two dense $q$-expansions
with $d$ terms. As $M_q$ is at least linear, the total complexity of
the series computations becomes
\[
O (\ell \, M_q (\ell d)).
\]

For functions $f$ related to $\eta$ quotients, such as $f = \wl^{2s}$ or
$f = j (z/\ell) = \frac {\left( \wt^{24} (z/\ell) + 16 \right)^3}{\wt^{24}
(z/\ell)}$, the effort for computing their $q$-expansions is negligible, since
it corresponds to a constant number of arithmetic operations with series
starting from the easily written down series expansion of $\eta$.

To write the $c_r$ as polynomials in $j$, one has to compute the non-positive
parts of the $q$-expansions of the $j^k$ for $1\leq k \leq d$. This can be
done in time $O (d \, M_q (d))$, which is negligible compared to the previous
steps because $d \in O (\ell)$. Identifying the $c_r$ as polynomials then
corresponds to solving triangular systems of linear equations, requiring
altogether $O \left( \ell d^2 \right)$ operations with integers.

The total complexity thus becomes
\[
O \left( \ell \, M_q (\ell^2) \right)
\]
integer operations.

Let $n \in O (\ell \log \ell)$ as in Section~\ref {ssec:height} be a bound
on the height of the modular polynomial. Then the bit complexity of the
algorithm is given by
\[
O \left( \ell^3 \log \ell \, M (n) \right)
\subseteq O \left( \ell^4 \log^{3 + \varepsilon} \ell \right)
\]

In \cite [Appendix A]{cl05} it is argued that in fact the running time is
higher, on grounds that the coefficients of the powers of $j$ grow faster than
$O (\ell \log \ell)$. However, this objection can be dismissed by carrying out
all computations modulo a sufficiently large prime of bit size
$n \in O (\ell \log \ell)$.
Alternatively, and probably preferably in practice, one may work modulo small
primes and use the Chinese remainder theorem.
Both approaches yield the desired complexity.

\begin {remark}
If $\Phi_\ell$ is sought only modulo some prime $p$, then the complete
algorithm can be carried out modulo $p$ as soon as $p$ is larger than $\ell$
to make the divisions in~\eqref {eq:newton} possible. The bit
complexity becomes
\[
O \left( \ell^3 \log \ell \log p \, (\log \log p)^{1 + \varepsilon} \right).
\]
\end {remark}

\subsection {Charles--Lauter}

In \cite {cl05}, the authors describe an algorithm to compute the classical
modular polynomial $\Phi_\ell$ directly modulo a prime $p$ without recourse to
$q$-expansions of modular functions. Instead, they rely on the moduli
interpretation of $\Phi_\ell$, which characterises pairs of $\ell$-isogenous
elliptic curves. So the algorithm manipulates only elliptic curves and
explicit isogenies over extension fields of $\F_p$.

Its basic building block is the computation of the instantiated polynomial
$\overline \Phi = \Phi_\ell (X, \jbar)$ with $\jbar \in \F_{p^2}$ the
$j$-invariant of a supersingular elliptic curve $E$. The roots of $\overline
\Phi$ are the $j$-invariants of the $\ell + 1$ curves that are
$\ell$-isogenous to $E$, and that may be obtained via Vélu's formul{\ae}
\cite {vel71} once the complete $\ell$-torsion of $E$ is known. These are
$\ell^2$ points defined over an extension of $\F_{p^2}$ of degree
$O (\ell)$ in the supersingular case; in fact, one expects a degree of $\Theta
(\ell)$ virtually all of the time. Thus writing down the $\ell$-torsion points
requires a time of $\Omega (\ell^3 \log p)$.

After having repeated the procedure for $\deg_j \Phi_\ell = \Theta (\ell)$
different values of $\jbar$, one obtains the coefficients of $\Phi_\ell$
modulo $p$ by interpolation. (This is analogous to the floating point approach
of Section~\ref {sec:evaluation}, and more details are given there.) The
complexity of the algorithm is thus at least
\[
\Omega \left( \ell^4 \log p \right),
\]
and an upper bound of
$O \left( \ell^{4+\varepsilon} \log^{2 + \varepsilon} p
+ \log^{4 + \varepsilon} p \right)$
is proved in \cite [Theorem~3.2]{cl05} under the generalised Riemann
hypothesis.

Hence, this algorithm is slower by a factor of order $\ell$ than the one
presented in Section~\ref {ssec:qexpansions} relying on $q$-expansions.
This is due to the costly determination of isogenies via their kernels: While
the isogenies are defined over $\F_{p^2}$, their kernels lie in an extension
of degree $O (\ell)$. The $q$-expansions of the conjugates of $j (z / \ell)$,
on the other hand, provide a synthetic description of the curves that are
$\ell$-isogenous to the one with $j$-invariant $j (z)$, and using them it is
sufficient to carry out all computations in $\F_p$ in order to obtain
$\Phi_\ell \bmod p$.

\section{Evaluation--interpolation}
\label{sec:evaluation}

The approach for calculating modular polynomials that is described in this
section is in fact neither new nor particularly involved. It has been
suggested to me by R.~Schertz during our work on the class invariants of \cite
{es04}, and later R.~Dupont pointed out to me that it can already be found in
\cite [Chapter~4.5]{bb87}. However, it has not been noticed before that the
algorithm allows to lower the exponent of the computational complexity by~$1$
and thus to reach an essentially optimal complexity if fast polynomial
arithmetic and fast techniques for evaluating modular functions as
described in \cite {dup07,eng06} are used.

The basic idea of the evaluation--interpolation approach is to specialise and
to compute the identity \eqref {eq:haupt} between modular functions
in several complex floating point arguments (the evaluation phase); and then
to interpolate the coefficients in order to recognise them as polynomials in
$j$. Precisely, \eqref {eq:haupt} can be rewritten as
\begin {equation}
\label {eq:special}
\Phi (X, j (z)) =
\prod_{\nu = 1}^{[\Gamma : \Gamma']} (X - f (M_\nu z))
= X^{[\Gamma : \Gamma']} + \sum_{r=1}^{[\Gamma : \Gamma']}
c_r (z) X^{[\Gamma : \Gamma'] - r}
\end {equation}
for all $z$ in the upper complex half plane.

Evaluating the conjugates $f (M_\nu z)$ of $f$ in a number of complex
arguments  $z_k$
with $\Im z_k > 0$, multiplying out the left hand side as a polynomial in $\C
[X]$ and separating the coefficients according to powers of $X$ yields the
values $c_r (z_k)$. Writing
$c_r = \sum_{s = 0}^{\deg_j \Phi} c_{r, s} j^{\deg_j \Phi - s}$,
the $c_r (z_k)$ are actually the values of the polynomial $c_r$ in $j (z_k)$,
so that the coefficients $c_{r, s} \in \C$ can be retrieved by interpolation
as soon as the $c_r (z_k)$ are known for $\deg_j \Phi + 1$ values of $z_k$.
The final step is to round the $c_{r, s}$ to rational integers. The degree of
the modular polynomial in $j$ or an upper bound on it may be known beforehand;
if it is not, then the algorithm can be repeated with increasing guesses for
$\deg_j \Phi$ until rounding to integers succeeds (and doubling the guess at
each time results in the same asymptotic complexity as taking the correct
value from the beginning).

Let $E (n)$ be the bit complexity for evaluating the modular functions $f$ or
$j$ at a precision of $n$ bits (here, $f$ is considered to be fixed, while $n$
tends to infinity with~$\ell$). Then the evaluation phase requires
$(\ell + 1) (\deg_j \Phi + 1)$ evaluations of $f$ at a cost of
\[
O \left( \ell \deg_j \Phi \, E (n) \right)
\]
and the reconstruction of $\deg_j \Phi + 1$ polynomials of degree $\ell + 1$
from their roots. Multiplying complex polynomials by the FFT, this
reconstruction step takes
\[
O \left( \deg_j \Phi \, \ell \log^2 \ell \, M (n) + \log n \, M (n) \right),
\]
where the (eventually negligible) term $\log n \, M (n)$ accounts for the
computation of a primitive root of unity of sufficiently high order to carry
out all the FFTs, and as usual $M (n)$ is the time needed to multiply two
$n$-bit numbers.

The interpolation phase consists of $\ell$ interpolations of polynomials of
degree $\deg_j \Phi$. Employing again fast algorithms such as \cite
[Algorithm~10.11]{gg99}, this takes
\[
O \left( \ell \deg_j \Phi \log^2 (\deg_j \Phi) \, M (n) \right),
\]
once the roots of unity are available.

It is shown in \cite {dup07} that among others, Dedekind's $\eta$ function
can be evaluated at precision $n$ in time
$E (n) \in O \left( \log n \, M (n) \right)$ uniformly for the argument in
$\Fc = \left\{ z \in \HH : |\Re z | \leq \frac {1}{2}, |z| \geq 1 \right\}$.
(The restriction to $\Fc$ is not crucial since we may freely choose
our interpolation points.)
So in particular, all functions built from $\eta$ such as $j$ and $\wl$ satisfy
$E (n) \in O \left( \log n \, M (n) \right)$. (Alternatively, one may use
multievaluation as described in \cite [Section~6.3]{eng06} with a slightly
worse amortised cost of $O \left( \log^2 n \, M (n) \right)$ per value,
which still allows to reach the final complexity of Result~\ref {th}.)

In total, the steps add up to a running time of
\begin {eqnarray*}
&& O \left( (\ell \deg_j \Phi \log^2 \max (\ell, \deg_j \Phi) + \log n) M (n)
\right) \\
& \subseteq &
O \left( (\deg_X \Phi \deg_j \Phi \log^2 \max (\deg_X \Phi, \deg_j \Phi) +
\log n) M (n) \right),
\end {eqnarray*}
and the latter formulation is valid for arbitrary congruence subgroups
$\Gamma'$ in the place of $\Gamma^0 (\ell)$.

Here, $n$ must be at least as large as the logarithmic height of $\Phi$, and
maybe bigger to account for rounding errors. In the case of the classical
modular polynomial $\Phi_\ell$ the bound $O (\ell \log \ell)$ of Section~\ref
{ssec:height} yields a complexity of
\[
O \left( \ell^2 \log^2 \ell M (\ell \log \ell) \right).
\]
This proves Result~\ref {th}.

\section{Further kinds of modular polynomials}
\label{sec:further}

\subsection {Modular functions of composite level}
\label {ssec:composite}

Besides its better complexity compared to the approach using $q$-expansions,
the eval\-u\-a\-tion--in\-ter\-po\-la\-tion algorithm has the advantage of a
great flexibility, which makes it easily adaptable to a variety of different
modular polynomials. In fact, the proof of Result~\ref {th} in Section~\ref
{sec:evaluation} does not use special properties of $\Gamma^0 (\ell)$ and is
valid for any congruence subgroup.

An application is given by the computation of the polynomial relationship
between $j$ and a modular function $f$ of composite level $N$. At first sight,
these polynomials are not of great interest. In the point counting or more
generally isogeny computation context, it is more efficient to express an
isogeny of composite degree as a composition of prime degree isogenies. But
this kind of modular polynomials occurs naturally in the context of \cite
{es04}, in which elliptic curves with complex multiplication are obtained via
modular functions $f$ of composite level. Levels $N = p q$ or $N = p^2$
with $p$ and $q$ prime are examined in \cite {es04}, but more general
settings can be
devised in a straightforward way. The polynomial relationship between $f$ and
$j$ is used to derive from a special value of $f$ the $j$-invariant of the
corresponding elliptic curve.

If the $q$-expansion approach were to be pursued to compute this kind of
modular polynomials, it would be necessary to somehow obtain the
$q$-expansions of the conjugates $f (M_\nu z)$ for a system of representatives
$(M_\nu)$ of $\Gamma' \backslash \Gamma$. This is straightforward only for
translations; for all other matrices, it is necessary to take the
transformation behaviour of $f$ under unimodular matrices into account, which
requires ad hoc computations for each particular function. (This is
illustrated by $j (S z / \ell) = j (\ell z)$ and $\wl^{2s} (Sz) = \ell^s
\left( \frac {\eta (\ell z)}{\eta (z)} \right)^{2s}$, which are not derived
from $j (z / \ell)$ resp. $\wl (z)$ by the same generic transformation.)
In the case of $\Gamma^0 (\ell)$, only the matrix $S$ asks for special
treatment; in the case of composite level, many more special matrices appear.
For instance, the full article \cite {es05} is concerned with deriving the
conjugates of only the functions
$\left( \frac {\eta (z / p_1) \eta (z / p_2)}{\eta (z / (p_1 p_2)) \eta (z)}
\right)^s$ for $p_1$, $p_2$ prime and
$s = 24 / \gcd (24, (p_1 - 1)(p_2 - 1))$.

In the evaluation--interpolation approach, the only difference between
$\Gamma^0 (\ell)$ and other congruence subgroups $\Gamma'$ is the enumeration
of the system of representatives $(M_\nu)$. For most interesting $\Gamma'$
(and in particular for arbitrary $\Gamma^0 (N)$), this is easy; in any case,
this step is independent of the particular function $f$. Then Result~\ref
{th} applies, and one obtains an algorithm that is essentially linear in its
output size.

\subsection {Schl\"afli equations}
\label {ssec:schlaefli}

In \cite {sch70}, Schl\"afli examines transformations of prime level~$\ell$ of
special modular functions different from $j$ that lead to particularly simple
modular equations. Weber gives a systematic treatment of them in \cite[\S\S
73--74]{web08}. Let $\f = \zeta_{48}^{-1} \frac {\eta ((z+1)/2)}{\eta (z)} =
\zeta_{48} \wt (z+1)$ be one of ``the Weber functions'', a modular function
of level $48$. Let $\ell$ be a prime not dividing this level. Then $g (z) = \f
(z/\ell)$ is the root of a monic polynomial $\Phi_\ell^\f (X)$ of degree $\ell +
1$ with coefficients in $\Z [\f]$.

The evaluation--interpolation approach allows to easily obtain this polynomial
with only minimal modifications to the algorithm. The function $\f$ being
modular for a subgroup of $\Gamma (48)$ and $\ell$ being different from $2$
and $3$, a quick computation reveals that $g$ is modular for $\Gamma (48) \cap
\Gamma^0 (\ell)$. The polynomial $\Phi_\ell^\f$ is thus given by an equation
analogous to~\eqref {eq:haupt}:
\begin {equation}
\label {eq:schlaefli}
\Phi_\ell^\f (X) = \prod_{\nu = 1}^{\ell + 1} (X - g (M_\nu z))
= X^{\ell + 1} + \sum_{r=1}^\ell c_r X^{\ell + 1 - r},
\end {equation}
where the $M_\nu$ range over a set of representatives of
$(\Gamma (48) \cap \Gamma^0 (\ell)) \backslash \Gamma (48)$. Such a set may be
obtained by multiplying the standard representatives of $\Gamma^0 (\ell)
\backslash \Gamma$ from the left by a matrix in $\Gamma^0 (\ell)$ such that
they end up in $\Gamma (48)$. Precisely, a possible set of representatives is
given by
$
\begin {pmatrix} 1 & 48 \nu \\ 0 & 1 \end {pmatrix}
$
for $\nu = 0, \ldots, \ell - 1$ (corresponding to the translations)
and
$
\begin {pmatrix} 1 - 48 k & 48 k \\ -48 k & 1 + 48 k \end {pmatrix}
$
with $k = 48^{-1} \bmod \ell$, corresponding to the inversion $S$.

So the only modification required for the evaluation--interpolation algorithm
to work is this adaptation of the matrices $M_\nu$, and of course $j$ has
to be
replaced by~$\f$ in the interpolation phase to recover the coefficients $c_r$
as elements of $\Z [\f]$.

Hence Result~\ref {th} clearly applies also to the Schl\"afli
equations, showing that they can be computed in essentially linear time with
respect to the output size.

Taking specifics of the function $\f$ into account, a more efficient algorithm
can be obtained, gaining a constant factor. Weber shows in
\cite [p.~266]{web08} that only every $24$-th coefficient in $\f$ of the $c_r$
may be non-zero. Precisely, $\f^i g^k$ having a non-zero coefficient implies
$\ell i + k \equiv \ell + 1 \pmod {24}$. Hence the number
of evaluations can be reduced by a factor of about $24$. (In \cite
[Ch.~4.5]{bb87}, similar sparse modular equations are studied, in which one
out of eight coefficients is non-zero. The Borwein's suggest in this chapter
the evaluation--interpolation approach while profiting of this sparseness.)
Weber shows that the $\Phi_\ell^\f$ are symmetric, which could also be taken
into account during the interpolation phase.

Weber derives the Schl\"afli equations in a very compact form. He considers
the new functions
\[
A = \left( \frac {\f}{g} \right)^r + \left( \frac {g}{\f} \right)^r
\text { and }
B = (\f g)^s + \legendre {2}{\ell} \frac {2^s}{(\f g)^s}
\]
with $r$, $s$ such that $12 | (\ell - 1) r$, $12 | (\ell + 1) s$ and
$\frac {(\ell - 1) r}{12} \equiv \frac {(\ell + 1) s}{12} \pmod 2$. Depending
on $\ell \bmod 24$, these conditions are satisfied by some $r$ and $s$ such
that $2r | \ell + 1$ and $2s | \ell - 1$. Weber \cite [p.~268]{web08} then
shows that
\[
\Phi_\ell^\f = (\f g)^{(\ell + 1) / 2} \left( A^{(\ell + 1) / (2r)}
- B^{(\ell - 1) / (2s)}
+ \sum_{\alpha = 0}^{(\ell + 1) / (2r) - 1} \;
\sum_{\beta = 0}^{(\ell - 1) / (2s) - 1} c_{\alpha, \beta} A^\alpha B^\beta
\right).
\]
In this form, the polynomial has about $\ell^2 / (4rs)$ coefficients, which is
often less than the roughly $\ell^2 / 24$ coefficients if it is written as an
equation between $\f$ and $g$. However, the more compact form appears to be
less suited to the evaluation--interpolation algorithm: There is no easy way,
in the spirit of \eqref {eq:haupt}, of obtaining its values when
interpreting it as a univariate polynomial in $A (z)$ with coefficients in $\Z
[B (z)]$. So instead of performing $O (\ell)$ interpolations of polynomials of
degree $O (\ell)$, one would apparently have to obtain all coefficients at once
by solving a linear system with $O (\ell^2)$ unknowns, resulting in a
complexity that is worse than that of Result~\ref{th}.

\subsection {Generalised Schl\"afli equations}
\label {ssec:general schlaefli}

The Schl\"afli equations of the previous paragraph can obviously be generalised
to further modular functions: Given a modular function $f$ and a prime $\ell$,
one may want to compute the polynomial relating $f (z)$ and $f (z/\ell)$.

In \cite[Ch.~5]{har04}, a theory analogous to Weber's is developped
for the functions $\wl$.
The derived polynomials are used to obtain symbolically minimal polynomials of
special algebraic values of $\eta$ quotients.

In the context of a $p$-adic algorithm for the computation of class polynomials
and ultimately elliptic curves with complex multiplications, another application
is developed in \cite [Ch.~6.8]{bro06}. Let $\tau$ be a quadratic integer. Then
the singular value $j (\tau)$ lies in the ring class field for the order $\Oc =
[1, \tau]_\Z$ and is the $j$-invariant of an elliptic curve with complex
multiplication by $\Oc$. In \cite {ch02}, the authors describe an algorithm to
compute the canonical lift of $j (\tau)$ to the $p$-adic numbers at arbitrary
precision by some kind of Newton iteration. The main algorithmic step consists
of applying an isogeny to $j (\tau)$ (or more precisely to the elliptic curve
with this $j$-invariant) that corresponds to a smooth principal ideal and that
can thus be realised by composing isogenies of manageable prime degrees. In
\cite {bro06}, the algorithm is generalised
to other class invariants. Let $f$ be a
modular function for $\Gamma^0 (N)$ such that $f (\tau)$ is a class invariant in
the sense that it also lies in the ring class field of $\Oc$. Let $\ell$ be a
prime number not dividing $N$ that splits as $(\ell) = \lf \overline \lf$ in
$\Oc$. As in the case of $j (\tau)$, one needs to ``apply an isogeny'' and
compute $f (\lf^{-1})$. This value is given as a root of the modular polynomial
$\Psi$ between $f (z)$ and $f (z/\ell)$, specialised with $f (\tau)$ in the
place of $f (z)$. (The modular polynomial $\Phi$ between $f (z)$ and $j (z)$,
treated in Section~\ref {ssec:composite},
also plays a role as it helps to choose whenever $\Psi$ has multiple roots: $f
(\lf^{-1})$ is also a root of $\Phi$ specialised with $j (\lf^{-1})$ in the
place of $j (z)$.)

A quick computation shows that $g (z) = f (z/\ell)$ is modular for $\Gamma^0
(\ell) \cap \Gamma^0 (N) = \Gamma^0 (\ell N)$, so that (\ref {eq:schlaefli})
yields the minimal polynomial of $g$ over $\C_{\Gamma^0 (N)}$ if the $M_\nu$ are
chosen as a system of representatives of $\Gamma^0 (\ell N) \backslash \Gamma^0
(N)$. Such a system is obtained precisely as for the Schl\"afli equations by
multiplying the standard representatives of $\Gamma^0 (\ell) \backslash \Gamma$
by matrices in $\Gamma^0 (\ell)$ so that they end up in $\Gamma^0 (N)$. For
instance, a set of representatives is given by
$
\begin {pmatrix} 1 & N \nu \\ 0 & 1 \end {pmatrix}
$
for $\nu = 0, \ldots, \ell - 1$
and
$
\begin {pmatrix} 1 - N k & N k \\ -N k & 1 + N k \end {pmatrix}
$
with $k = N^{-1} \bmod \ell$.
So at first sight, the evaluation--interpolation algorithm seems to apply.

The problem lies in recognising the coefficients $c_r$ as polynomials in $f$. In
fact, the $c_r$ are modular for $\Gamma^0 (N)$, so unless $C_{\Gamma^0 (N)} = \C
(f)$ (otherwise said, $f$ is a \textit {Hauptmodul} for $\Gamma^0 (N)$),
there is no reason that they are rational in $f$. Generically, one expects
$C_{\Gamma^0 (N)} = \C (j, f)$ and may hope (under suitable integrality and
rationality conditions) to recover the $c_r$ as elements of $\Z [j, f]$.

A necessary condition for $f$ being a Hauptmodul is that the modular curve $X_0
(N)$ is of genus~$0$, which limits the possible values of $N$ to a very short
list. In general, it will be necessary to replace the minimal polynomial of $g$
over $\C_{\Gamma^0 (N)}$ by that over $\C (f)$; its degree will be $\ell + 1$
multiplied by $[\C_{\Gamma^0 (N)} : \C (f)]$. It is not clear whether the
evaluation--interpolation approach can be adapted to this context, since $\C
(f)$ is not the field of modular functions for any congruence subgroup.

\section{Implementation}
\label{sec:implementation}

The evaluation--interpolation algorithm has been implemented in C by the author
for different kinds of modular polynomials. One big advantage of the algorithm
(besides its superior complexity) is that it is rather straightforward to
implement, especially in the context of complex multiplication (cf. \cite
{eng06}), where the major building blocks such as the evaluation of modular
functions and computations with polynomials over floating point numbers are
already in use. The author's implementation relies on the existing libraries
\texttt {gmp} \cite {gmp} with an assembly patch for AMD64 by P.~Gaudry \cite
{gau05ass} for fast basic multiprecision arithmetic, \texttt {mpfr} \cite
{mpfr} and \texttt {mpc} \cite {mpc}, which provide elementary operations
with arbitrary precision floating point real and complex numbers with exact
rounding. A library for fast polynomial operations via Karatsuba, Toom--Cook
and the FFT has been written for the occasion \cite{mpfrcx}.
All timings mentioned in the
following have been obtained on Opteron processors clocked at 2.4~GHz. The
heights are given as logarithms in base $2$, and the arithmetic precision is
given in bits.

\subsection{Details of the implementation}

\subsubsection{Details slowing the implementation down asymptotically}

Of the different techniques for evaluating modular functions described in \cite
{eng06}, the asymptotically fast ones are not beneficial in the present context.
For computing univariate class polynomials, they do not pay off at degrees
below 100~000. Such high degrees are for the time unachievable for
bivariate modular polynomials, since the number of coefficients would make it
impossible to store the result. Hence all examples have been computed via the
sparse representation of the $\eta$ function as a $q$-series.

The interpolation phase turns out to take negligible time; thus a simple
quadratic algorithm using iterated differences is employed instead of a
quasilinear one.

\subsubsection{Details speeding the implementation up}

In order to work with real instead of complex numbers as much as possible,
purely imaginary values are chosen for the $z_k$. Then the $q (z_k)$ are real,
and so are the values $j (z_k)$ and thus, by (\ref {eq:special}), the
$\Phi (X, j (z_k))$. While the values of the conjugates $f (M_\nu z_k)$ still
have to be computed as complex numbers, roughly half of them suffice for
functions $f$ for $\Gamma^0 (\ell)$ with rational $q$-coefficients: $f (z_k)$ is
real, $f (z_k - r)$ is the complex conjugate of $f (z_k + r)$ for $r \in
\Z$, and finally $f (-1/z_k)$ has to be real again to obtain a real polynomial
in the end. When reconstructing $\Phi (X, j (z_k))$ from its roots, the complex
conjugate roots can be grouped so that only real arithmetic is required.

Also the interpolation phase uses only real numbers. It is furthermore helped by
choosing the $z_k$ such that the $j (z_k)$ lie in a simple arithmetic
progression, for instance, $j (z_k) = 1728 + k$ for $k \geq 1$.

\subsubsection{Parallelisation}

Another nice feature of the evaluation--interpolation algorithm is that it can
be parallelised almost arbitrarily, in principal up to distributing at the level
of each modular function evaluation. A more natural, coarser parallelisation
of the evaluation phase, which needs almost no communication, is amply
sufficient in practice: Each processor treats a different $z_k$ and computes an
approximation $\Phi (X, j (z_k))$ by evaluating the conjugates $(f (M_\nu
z_k))_\nu$ and reconstructing the polynomial from these roots. The result is
written into a file on a shared hard disk. The interpolation phase has been
too fast to warrant a parallel implementation. If it were to become a
bottleneck, a natural way of distributing the work would be to have each
processor compute the coefficients of a different polynomial $c_r \in \Z [j]$.

Implementing the communication via files solves a second problem. During the
evaluation phase, the matrix $\left( c_r (z_k) \right)_{k, r}$ is computed
row-wise; for the interpolation, it is accessed column-wise. However, the matrix
is roughly as big as the final modular polynomial and for the largest computed
examples (see Section~\ref {ssec:examples0}) does not fit into main memory any
more. Writing each row into a different file during the evaluation and reading
all files in parallel during the interpolation can be seen as a means of
transposing the matrix on disk without having to keep it in main memory.

\subsubsection{Arithmetic precision}

Except for the classical modular polynomials, no upper bound on the height of
the polynomials is known, and even in the classical case, the bound is not
completely explicit, but contains an $O (1)$ term (see Section~\ref
{ssec:height}). In the implementation, an approximate bound is obtained by
carrying out all computations first at very low precision ($100$ bits). Due to
the numeric instability of interpolation, the resulting coefficients are wrong
already in the first digit, but one obtains an idea of the height of the correct
polynomial (actually, the height at low precision turns out to be larger than
the real one). This estimate is multiplied by a factor (between $1.1$ and $2$
depending on the function and the level), determined experimentally, to deal
with rounding errors.

In the context of class polynomial computations, it has been observed in \cite
{eng06} that virtually no rounding errors occur: Evaluating modular functions
via the $q$-development of $\eta$ as well as multiplication of polynomials is
rather uncritical from a numerical point of view. It suffices to add a few bits
to the target precision. This should also hold for the evaluation phase in our
context. In practice, it can be observed, however, that a significantly higher
precision is needed than just the height of the result. This can only be due to
the interpolation phase, which implicitly handles Vandermonde matrices, that are
known to be ill conditioned. It is an open problem to choose the interpolation
points $j (z_k)$ so as to minimise the rounding errors.

\subsection{Examples}
\label {sec:examples}

\subsubsection{Modular polynomials for $\Gamma^0 (\ell)$}
\label {ssec:examples0}

For achieving the point counting record for elliptic curves with the
Schoof--Elkies--Atkin algorithm \cite {em06}, the polynomials have been
systematically computed with the function
\begin {equation}
\label {eq:atkin}
f_{\ell, r} = \frac {T_r ( \eta (z) \eta (\ell z))}{\eta (z) \eta (\ell z)}
\end {equation}
(evaluated in $-1/z$) suggested in \cite {mul95alg}. Here,
\[
T_r (f) = \frac {1}{r} \sum_{\nu = 0}^{r - 1} f \left( \frac {z + 24 \nu}{r}
\right) + f (rz)
\]
is the $r$-th Hecke operator for a prime $r \geq 5$ such that
$24 \divides (r - 1)(\ell + 1)$, $r$ is a quadratic residue modulo $\ell$ and
$\ell$ is a quadratic residue modulo $r$. (The unusual factor $24$ in front of
$\nu$ takes care of the fact that the exponents in the $q$-expansion of $\eta$
are not integral due to the factor $q^{1/24}$; it eliminates $24$-th roots of
unity.)

It is not advisable to obtain the values of $f_{\ell, r}$ via its
$q$-expansion. First of all, the expansion is dense and thus much slower to
evaluate than that of $\eta$. But more importantly, $q$ tends to infinity when
$z$ approaches the cusp $0$ of the fundamental domain for $\Gamma^0 (\ell)$,
so that the $q$-series converges arbitrarily slowly. This is not an issue when
relying on the evaluation of $\eta$ only, since its known transformation
behaviour under unimodular matrices allows to transform the arguments into the
fundamental domain for $\Gamma$, whose only cusp is at infinity and corresponds
to $q = 0$. So the implementation evaluates (\ref {eq:atkin}) numerically with
$2 r + 4$ evaluations of $\eta$ per value of $f_{\ell, r}$. The resulting
dependence of the running time on $r$ can be seen clearly in the following
examples.

\begin {center}
\begin {tabular}{lr}
$\ell$           & $2039$ \\
\hline
$r$              & $5$ \\
degree in $j$    & $136$ \\
\hline
estimated height & $5816$ \\
precision        & $6397$ \\
height           & $5040$ \\
\hline
time for evaluation    & $5900$ s \\ 
time for interpolation & $110$ s\\
\end {tabular}
\hspace {1cm}
\begin {tabular}{lr}
$\ell$           & $2017$ \\
\hline
$r$              & $61$ \\
degree in $j$    & $156$ \\
\hline
estimated height & $6211$ \\
precision        & $6832$ \\
height           & $5311$ \\
\hline
time for evaluation    & $74000$ s \\ 
time for interpolation & $130$ s\\
\end {tabular}
\end {center}

The largest computed polynomial has a level of $10079$. The computation takes
less than two weeks on a small cluster with ten processors; as a compressed
text file, the polynomial fills about $16$~Gb of space. This illustrates
that the limiting factor for the algorithm becomes not so much its running
time, but rather the space requirements of its output, as can be expected for
algorithms that have an essentially linear time complexity with respect to
their output size.

\begin {center}
\begin {tabular}{lr}
$\ell$           & $10079$ \\
\hline
$r$              & $5$ \\
degree in $j$    & $672$ \\
\hline
estimated height & $31865$ \\
precision        & $35051$ \\
height           & $28825$ \\
\hline
time for evaluation    & $10~000~000$ s $\approx 120$ d \\ 
time for interpolation & $56~000$ s $\approx 16$ h \\
\end {tabular}
\end {center}

It would be natural to try to linearly combine functions $f_{\ell, r}$ with
different $r$ to reduce the pole order at infinity, or even to find a function
on $X_0 (\ell)$ with minimal pole order. As a result, the degree of the modular
polynomial in $j$ (and thus the number of interpolation points) and the size of
its coefficients (and thus the required precision) could be reduced. A high
price, however, would have to be paid by an increased running time for the
evaluation. As high performance clusters are becoming increasingly available,
while bandwidth and storage space appear to be the bottleneck for modular
polynomial computation, this might, however, be the road to follow for
constructing a manageable database of modular polynomials up to a high level.

\subsubsection {Schl\"afli equations}

As explained in Section~\ref {ssec:schlaefli}, a trivial modification to the
enumeration of the $\ell + 1$ matrices $M_\nu$ adapts the code to computing the
Schl\"afli equations between $\f (z)$ and $\f (z / \ell)$. It is then
imperative to change the interpolation code so as to take advantage of the
sparseness of the polynomial, in which only every $24$-th coefficient is
non-zero. This reduces the number of evaluations and (after a suitable
transformation) the degree of the interpolated polynomials roughly by a factor
of $24$. The latter is also important since the interpolation phase becomes
numerically more stable.

\begin {center}
\begin {tabular}{lr}
$\ell$           & $2039$ \\
\hline
number of interpolation points    & $86$ \\
\hline
estimated height & $2829$ \\
precision        & $3111$ \\
height           & $2380$ \\
\hline
time for evaluation    & $230$ s \\
time for interpolation & $35$ s\\
\end {tabular}
\end {center}

\subsubsection {Generalised Schl\"afli equations}

As argued in Section~\ref {ssec:general schlaefli}, the minimal polynomial of $g
(z) = f (z / \ell)$ over $\C (f)$ for some function $f$ for $\Gamma^0 (N)$ need
not be of degree $\ell + 1$ any more, in which case it is a priori unclear how
to obtain it by evaluation--interpolation. But besides the functions $f$
that are Hauptmoduln for $\Gamma^0 (N)$, a few others may be handled by this
approach. Namely, let
\[
\w_{p_1, p_2} = \frac {\eta (z / p_1) \eta (z / p_2)}
{\eta (z) \eta (z / (p_1 p_2))}
\]
with $p_1$, $p_2$ prime and $24 \divides (p_1 - 1)(p_2 - 1)$ be the functions of
level $N = p_1 p_2$ suggested as class invariants in \cite {es04}. These
functions are invariant under the Fricke--Atkin--Lehner involution. Now, $X_0^+
(N)$ is of genus $0$ for $N = 35$ and $N = 39$ \cite {ogg74}, and apparently,
the double $\eta$ quotients are Hauptmoduln in this case. The degree of the
modular equation of transformation level $\ell$ thus remains $\ell + 1$ (an
observation made for $N = 35$ in \cite [Ch.~7.4]{bro06}).

For instance, with $f = \w_{3,13}$ and $g (z) = f (z / \ell)$, one obtains the
following polynomials $\Phi_\ell^{\w_{3,13}}$:
\begin {eqnarray*}
\Phi_2^{\w_{3,13}} & = &
g^3 + f^3 - g^2 f^2 - g f + 2 (g^2 f + f g^2) \\
\Phi_5^{\w_{3,13}} & = &
g^6 + f^6 - g^5 f^5 - g f
+ 35 g^3 f^3 \\
&& + 5 (g^5 f^4 + g^4 f^5 + g^5 f + g f^5 - g^4 f^3 - g^3 f^4
     - g^3 f^2 g^2 f^3 + g^2 f + g f^2) \\
&& + 10 (- g^5 f^2 - g^2 f^5 + g^4 f^4 + g^4 f^2 + g^2 f^4 - g^4 f - g f^4
      + g^2 f^2) \\
\Phi_7^{\w_{3,13}} & = &
g^8 + f^8 - g^7 f^7 -g f \\
&& + 7 (g^7 f^6 + g^6 f^7 - g^7 f^5 - g^5 f^7 + g^6 f^4 + g^4 f^6
        + g^6 f^2 + g^2 f^6) \\
&& + 7 (g^4 f^2 + g^2 f^4 - g^3 f - g f^3
        + g^2 f + g f^2) \\
&& + 14 (- g^7  f  - g f^7 + g^5 f^4 + g^4 f^5 + g^5 f^3 + g^3 f^5
         + g^4 f^3 + g^3 f^4) \\
&& + 21 (- g^7 f^4 - g^4 f^7 + g^7 f^3 + g^3 f^7 - g^6 f^5 - g^5 f^6
         + g^6 f^3 + g^3 f^6) \\
&& + 21 (g^5 f^2 + g^2 f^5 + g^5 f + g f^5
         - g^4 f - g f^4 - g^3 f^2 - g^2 f^3) \\
&& + 42 (g^6 f^6 + g^2 f^2) - 182 g^4 f^4
\end {eqnarray*}

Like the Schl\"afli equations, the polynomials are symmetric in $f$ and $g$,
but as can be seen from these and further examples, they are no more sparse.

\bibliographystyle {amsplain}
\bibliography {modcomp}
\end{document}